\documentclass[12pt,a4paper]{article}
\usepackage[top=2cm, bottom=2.5cm, left=1.6cm, right=1.6cm]{geometry}
\usepackage{amsmath,amssymb,amsfonts,amsxtra,graphics,graphicx,amsthm,bbm,array}
\usepackage{mathrsfs}
\usepackage{hyperref}
\usepackage{CJK}

\usepackage[T1]{fontenc}
\usepackage[utf8]{inputenc}
\usepackage{authblk}
\usepackage{dcolumn}
\usepackage[vcentermath]{youngtab}

\newcommand{\bC}{\ensuremath{\mathbb{C}}}

\newcommand{\bZ}{\ensuremath{\mathbb{Z}}}

\newcommand{\scM}{\ensuremath{\mathscr{M}}}

\newcommand{\fraksl}{\ensuremath{\mathfrak{sl}}}

\newcommand{\cN}{\mathcal{N}}

\newcommand{\cW}{\mathcal{W}}

\newcommand{\SL}{\mathrm{SL}}

\newcommand{\Li}{{\rm Li}}

\newcommand{\smalltet}[2]{\raisebox{-#2cm}{\includegraphics[width=#1cm]{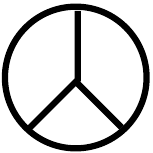}}}
\newcommand{\smalltrivalentA}[2]{\raisebox{-#2cm}{\includegraphics[width=#1cm]{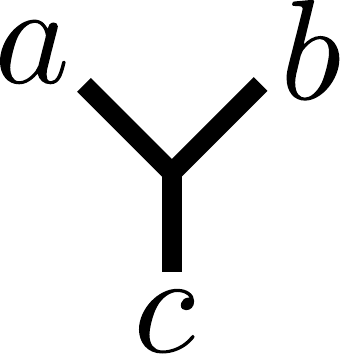}}}

\newcommand{\half}{\frac{1}{2}}

\def\beaa{\begin{eqnarray*}}
\def\eeaa{\end{eqnarray*}}
\def\bee{\begin{equation*}}
\def\eee{\end{equation*}}
\def\bea{\begin{eqnarray}}
\def\eea{\end{eqnarray}}
\def\be{\begin{equation}}
\def\ee{\end{equation}}
\def\ba{\begin{align}}
\def\ea{\end{align}}

\newcommand{\bem}{\begin{pmatrix}}
\newcommand{\eem}{\end{pmatrix}}

\def\={\;  = \;}
\def\+{\, + \,}

\def\wh{\widehat}

\def\rt2{\sqrt{2}}

\begin{document}
\title{\begin{flushright}\begin{footnotesize}
{NIKHEF-2014-008}\end{footnotesize}
\end{flushright}\bigskip
Trivalent graphs, volume conjectures and character varieties}

\author[1]{Satoshi Nawata}
\author[2]{Ramadevi Pichai}
\author[2]{Zodinmawia}
\affil[1]{NIKHEF theory group, Science Park 105, 1098 XG Amsterdam, The Netherlands}
\affil[2]{Department of Physics, Indian Institute of Technology Bombay,  India, 400076}
\date{}

  \maketitle

\abstract{The generalized volume conjecture and the AJ conjecture (a.k.a. the quantum volume conjecture) are extended to $U_q(\fraksl_2)$ colored quantum invariants of the theta and tetrahedron graph. The $\SL(2,\bC)$ character variety of the fundamental group of the complement of a trivalent graph with $E$ edges in $S^3$ is a Lagrangian subvariety of the Hitchin moduli space over the Riemann surface of genus $g=E/3+1$. For the theta and tetrahedron graph, we conjecture that the configuration of  the character variety is locally determined by large color asymptotics of the quantum invariants of the trivalent graph in terms of complex Fenchel-Nielsen coordinates. Moreover, the $q$-holonomic difference equation of the quantum invariants provides the quantization of the character variety. \\}

{\bf Mathematical Subject Classification (2010)}: 14D21, 53D30, 57M50, 81R50\\

{\bf Keywords}:  quantum invariants of trivalent graphs, generalized volume conjecture, AJ conjecture, character varieties


\vspace{.5cm}

\section{Introduction and summary}
The volume conjecture provides a relation between colored Jones 
polynomials of a knot and the hyperbolic volume of the complement of the knot in $S^3$. After its formulation \cite{Kashaev:1996kc,Murakami:1999}, we have witnessed its extensions of various kinds which relate quantum invariants of knots and classical geometry of knot complements. In particular, the generalized volume conjecture \cite{Gukov:2003na} connects quantum invariants and character varieties of knot groups, and the AJ conjecture (a.k.a. the quantum volume conjecture) \cite{Garoufalidis:2003a,Garoufalidis:2003b,Gukov:2003na} gives the quantization of character varieties. In this letter, we will consider $U_q(\fraksl_2)$ quantum invariants of simple planar graphs in the context of the generalized volume conjecture and the AJ conjecture.

A trivalent graph $\Gamma$  is a closed one-dimensional cell complex where three edges meet at each vertex. Therefore, the number $E$ of edges and that $V$ of vertices are related by $2E=3V$. Colored $U_q(\fraksl_2)$ quantum invariants $J_{n_1,\cdots,n_E}(\Gamma;q)$ of a trivalent graph \cite{Constantino:2009,masbaum1994,Roland} can be obtained as follows. First, each edge of the graph carries a representation of  $U_q(\fraksl_2)$ as a color which is specified by an integer $n_i\in \bZ$ corresponding to the spin-$n_i/2$ representation of $U_q(\fraksl_2)$. The colors ($a,b,c$) on  three edges meeting a trivalent vertex $\smalltrivalentA{.6}{.15}$ must obey the \emph{fusion rule}, \textit{i.e.} $|a-b|\leq c\leq a+b$. We also call a pair $(a,b,c)$ subject to the fusion rule an \emph{admissible set}. Then, the quantum invariant of the colored trivalent graph can be obtained by summing over angular momentum states of $U_q(\fraksl_2)$ Clebsch-Gordan ($q$-CG) coefficients \cite{Kirillov:1988} associated to trivalent vertices. Specifically, this letter will deal with the first two simplest trivalent graphs: a theta graph (Figure \ref{fig:theta2}) and a tetrahedron graph (Figure \ref{fig:tetrahedron}).

In fact, colored quantum invariants of trivalent graphs have  been taken into account in terms of the volume conjecture \cite{Constantino:2007,Constantino:2010,Constantino:2014,Murakami:2004,Murakami,Roland}. Namely, a suitable large color limit of quantum invariants of a trivalent graph can be related to the hyperbolic volume of the complement of the graph. Hence, our aim is to extend the generalized volume  conjecture and the AJ conjecture (the quantum volume conjecture)  for trivalent graphs.

For completeness, we will first briefly review the generalized volume conjecture and the AJ conjecture for colored Jones polynomials of knots. The generalized volume conjecture \cite{Gukov:2003na} amounts to the statement that   in the double scaling limit, $n\rightarrow\infty$, $q=e^{\hbar}\rightarrow 1$ with $x=q^{n/2}$ fixed, the asymptotic behavior of the 
colored Jones polynomials $J_n(K;q)$, for any knot $K$ takes the form
\be
J_n(K;q)\;\overset{{n \to \infty \atop \hbar \to 0}}{\sim}\;
\exp\left( \frac{1}{\hbar} \int \log y \frac{dx}{x} \,+\, \ldots \right)~,
\label{VC}
\ee
where the $y(x)$ essentially gives the zero locus of the classical $A$-polynomial $A(K;x,y(x))=0$. The $A$-polynomial is also a knot invariant which describes the moduli space of $\SL(2,{\bC})$ flat connections over the knot complement $S^3\backslash K$. 
This moduli space  is also called  as 
$\SL(2,{\bC})$ character variety of the fundamental group of the knot complement $S^3\backslash K$. Recall that the moduli space $\scM_{\rm flat}(T^2,\SL(2,{\bC}))$ of $\SL(2,{\bC})$ flat connections on the boundary torus of the knot is a hyper-K\"ahler variety $\bC^\times \times\bC^\times /\bZ_2$.  The space $\bC^\times  \times\bC^\times$ is spanned by the holonomy eigenvalues of the $\SL(2,{\bC})$ gauge connection along the meridian $x$ and the longitude $y$ of the torus, and $\bZ_2$ is the Weyl group symmetry of the gauge group $\SL(2,{\bC})$.  The moduli space of $\SL(2,{\bC})$ flat connections on the knot complement $S^3\backslash K$ is the Lagrangian subvariety of the moduli space $\scM_{\rm flat}(T^2,\SL(2,{\bC}))$, with respect to the symplectic form $\omega=\frac{1}{\hbar}d\log{x}\wedge d\log y$,  defined by the zero locus of the $A$-polynomial:
\bea
\scM_{\rm flat}(S^3\backslash K,\SL(2,{\bC}))=\{ (x,y)\in\scM_{\rm flat}(T^2,\SL(2,{\bC}))| A(K;x,y)=0\}~.
\eea
Furthermore, the quantization of character varieties of knot groups $\pi_1(S^3\backslash K)$ can be obtained by considering the recursion relations of Jones polynomials, which is called the AJ conjecture \cite{Garoufalidis:2003a,Garoufalidis:2003b} or the quantum volume conjecture \cite{Gukov:2003na}. The recursion of the colored Jones polynomials of minimal order can be manipulated in the form
\be
 {\widehat A}(K;\hat x,\hat y; q)J_n(K;q)=0~,
\ee
where the operator $\hat x$ and $\hat y$ has the following action on the 
colored Jones polynomial:
\bea\label{xyactionJ}
\hat x J_n(K;q)=q^{n/2} J_n(K;q)~, \qquad
\hat y J_n(K;q)=J_{n+1}(K;q)~.
\eea
The quantum $A$-polynomial $\wh A(K;\hat x,\hat y;q)$ reduces to the classical one at $q=1$: $\wh A(K; x, y;q=1)=A(K;x,y)$.

It is natural to ask whether the generalized volume conjecture and the AJ conjecture are
applicable for quantum invariants of a trivalent graph $\Gamma$ with $E$ edges colored by representations
$n_1,\cdots,n_E$. Although one can consider knotted trivalent graphs in general, we restrict ourselves to the case of planar graphs, graphs which can be projected to the plane, for brevity in this paper. For a planar trivalent graph, the boundary of the tubular neighborhood of the graph $\Gamma$ is no longer a torus but a Riemann surface of genus $g\ge2$. Actually, the genus is expressed by the number of edges via $g=E/3+1$. This can be explained as follows. Let us project the trivalent graph $\Gamma$ on two-sphere $S^2$. This provides a cell decomposition of $S^2$ where the number $F$ of two-complex cells  is indeed $F=g+1$. Since the Euler characteristics of $S^2$ is two, we have the relation 
\bea
V-E+F=2   ~.
\eea
Using $2E=3V$ and $F=g+1$, we obtain $g=E/3+1$. In fact, the moduli space $\scM_{\rm flat}(\Sigma_g,\SL(2,{\bC}))$  of $\SL(2,{\bC})$ flat connections on the boundary Riemann surface of genus $g$ is known as the Hitchin moduli space on $\Sigma_{g}$, which is a non-trivial hyper-K\"ahler variety of $\dim_\bC=6g-6=2E$ \cite{Hitchin:1987}. Assigning a simple closed curve corresponding to each edge of the graph to the Riemann surface $\Sigma_{g}$, pants decomposition of $\Sigma_{g}$ is obtained. (See for example, Figure \ref{fig:theta2} and \ref{fig:tetrahedron}.) As a result, the set of the $\SL(2,{\bC})$ holonomy eigenvalue $x_i$ along a simple closed curve and its canonical conjugate ``twist'' variable $y_i$ becomes local Darboux coordinates of the Hitchin moduli space, so-called complex Fenchel-Nielsen coordinates \cite{Goldman,Dimofte:2013lba,Dimofte:2014ria}. Since the moduli space $\scM_{\rm flat}(S^3\backslash \Gamma,\SL(2,{\bC}))$  of $\SL(2,{\bC})$ flat connections on the complement of the trivalent graph $\Gamma$ in $S^3$ is a Lagrangian subvariety of $\scM_{\rm flat}(\Sigma_g,\SL(2,{\bC}))$, the configuration of the moduli space $\scM_{\rm flat}(S^3\backslash \Gamma,\SL(2,{\bC}))\subset \scM_{\rm flat}(\Sigma_g,\SL(2,{\bC}))$ can be locally specified by complex Fenchel-Nielsen coordinates. Now, one can ask a question whether there is a relation between colored quantum invariants $J_{n_1,\cdots,n_E}(\Gamma;q)$ of a graph and the the configuration of the moduli space $\scM_{\rm flat}(S^3\backslash \Gamma,\SL(2,{\bC}))$.

In what follows, we will make a modest step to answer to this question for the theta and tetrahedron graph. Therefore, let $\Gamma$ be either the theta or tetrahedron graph and  $J_{n_1,\cdots,n_E}(\Gamma;q)$  be its colored quantum invariants.
In the double scaling limit $q=e^\hbar\to1$ and $n_i\to \infty$ with $q^{n_i/2}=x_i$ fixed $(i=1,\cdots,E)$, we conjecture that the large color asymptotics of quantum invariants of $\Gamma$ exhibits the form
\bea
 J_{n_1,\cdots,n_E}(\Gamma;q) &\overset{{n_i \to \infty \atop \hbar \to 0}}{\sim}& e^{\tfrac{1}{\hbar}W(\Gamma;\{x_i\})}
=\exp\left( \frac{1}{\hbar} \int\sum_i \log y_i \frac{dx_i}{x_i} \,+\, \ldots \right)~~. 
\label{vcj-graph}
\eea
where the integral is performed over the moduli space $\scM_{\rm flat}(S^3\backslash \Gamma,\SL(2,{\bC}))$   determined by $3g-3=E$ set of equations
\bea\label{character-variety}
A_j(\Gamma; y_j,x_1,\ldots, x_E)=0 ~, \quad j=1,\cdots, E~.
\eea
Note that $\{x_i,y_i\}_{i=1,\cdots,E}$ are complex Fenchel-Nielsen coordinates of the Hitchin moduli space. This set of equations can be obtained by solving the equations
\be
y_j=\exp\left(x_j\frac{\partial W(\Gamma;\{x_i\})}{\partial x_j}\right)~. 
\ee
Furthermore, the moduli space $\scM_{\rm flat}(S^3\backslash \Gamma,\SL(2,{\bC}))$ defined by \eqref{character-variety}  is the Lagrangian subvariety of $ \dim_\bC=3g-3=E$  in the moduli space $\scM_{\rm flat}(\Sigma_g,\SL(2,{\bC}))$  with respect to the symplectic form $\omega=\frac{1}{\hbar}\sum_{i=1}^E d\log{x_i}\wedge d\log y_i$.

It is straightforward to formulate the AJ conjecture for the theta and tetrahedron graph. 
Similar to the operator actions $\hat x$ and $\hat y$ in \eqref{xyactionJ}, let us define the acton of operators $\hat x_j$ 
and $\hat y_j$ ($j=1,...,E$) on the graph invariant as follows:
\bea\label{xyactionJ2}
\hat x_j J_{n_1,\cdots,n_E}(\Gamma;q)&=&q^{n_j/2} J_{n_1,\cdots,n_E}(\Gamma;q) \cr
\hat y_j J_{n_1,\cdots,n_E}(\Gamma;q)&=&J_{n_1,\cdots,n_j+2,\cdots,n_E}(\Gamma;q),
\eea
We should stress out that the fusion rule allows to increase the rank of a representation by two boxes in terms of Young tableaux unlike the case of knot invariants. Then, a recursion relation of minimal order with respect to one of colors  
 \be
b_k^{(j)}(\hat x_1,...,\hat x_E,q)J_{n_1,\cdots,n_j+2k,\cdots,n_E}(\Gamma;q)+\cdots+b_0^{(j)}(\hat x_1,...,\hat x_E,q)J_{n_1,\cdots,n_j,\cdots,n_E}(\Gamma;q)=0~,\label{q-holo}
\ee
can be written as
\be
 {\widehat A}_j(\Gamma;\hat y_j,\hat x_1,\ldots,\hat x_E; q)J_{n_1,\cdots,n_E}(\Gamma;q)=0~,\label{gamaaj}
\ee
where 
\be
 {\widehat A}_j(\Gamma;\hat y_j,\hat x_1,\ldots,\hat x_E ;q)=\sum_{\ell=0}^k b_\ell^{(j)}(\hat x_1,...,\hat x_E,q)\hat y_j^\ell~.\label{gamaaj}
\ee
The set of the operators $ {\widehat A}_j(\Gamma;\hat y_j,\{\hat x_i\}; q)$  ($i=1,...,E$)  provides the quantization of the character variety $\scM_{\rm flat}(S^3\backslash \Gamma,\SL(2,{\bC}))$ defined by \eqref{character-variety} in the sense that 
\be\label{classical-limit}
 {\widehat A}_j(\Gamma;y_j,\{x_i\}; q=1)=A_j(\Gamma; y_j,\{x_i\})~.
\ee
In the following sections, we will demonstrate explicit calculations of $A_j(\Gamma; y_j,\{x_i\})$ for the theta and tetrahedron graph.

\section{Theta graph}
\begin{figure}[h]
\begin{center}
\includegraphics[scale=.35]{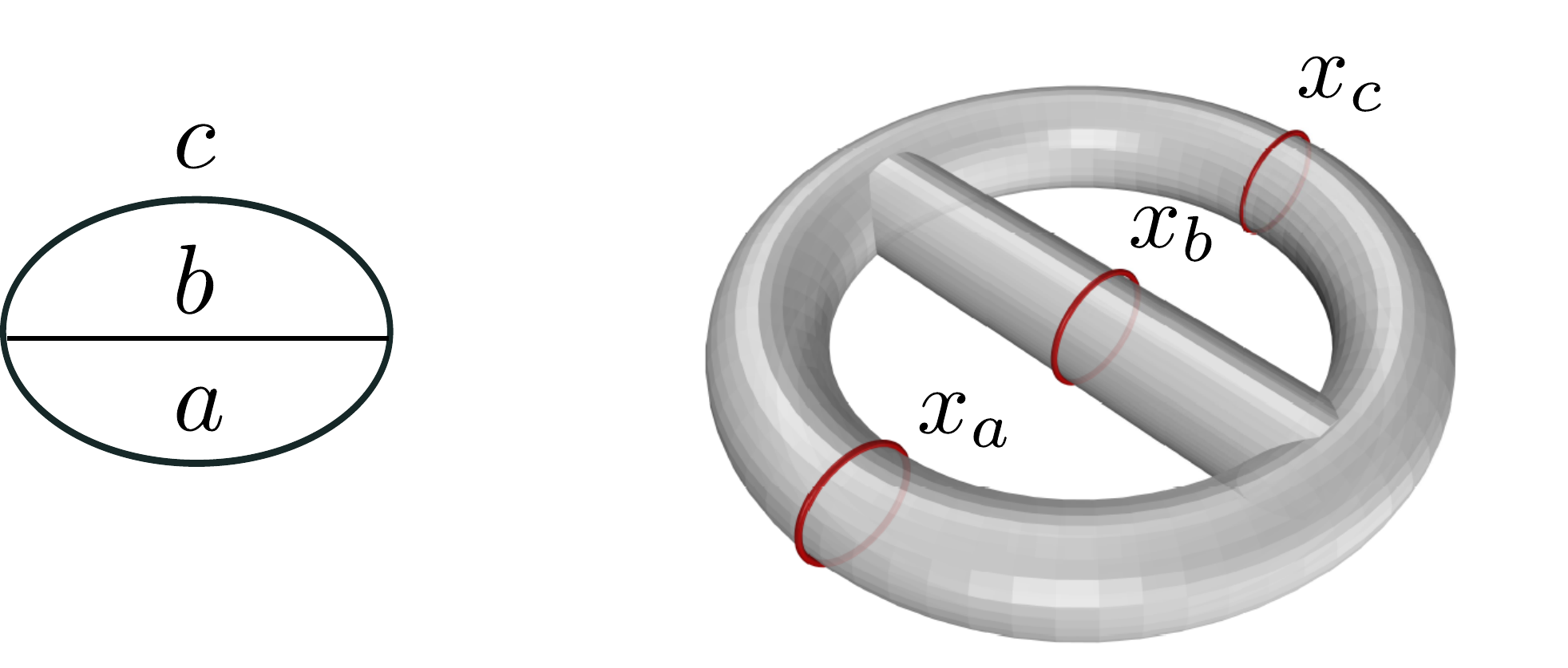}
\end{center}
\caption{Theta graph and its tubular neighborhood.}
\label{fig:theta2}
\end{figure}

 A theta graph is  the simplest trivalent graph with  three edges $E=3$  and two trivalent vertices $V=2$ as shown in Figure \ref{fig:theta2}.   The colored quantum invariant $J_{a,b,c}(\Theta;q)$ ($a,b,c\in\bZ$) of the theta graph 
can be written as a state sum involving product of two Clebsch-Gordan coefficients \cite{Constantino:2009}. In \cite{masbaum1994}, the invariants of the theta graph is written in the succinct form by using the Wenzl's recursion formula for Jones idempotents of the Temperly-Lieb algebra
\bea\label{theta-inv}
J_{a,b,c}(\Theta;q) &=&(-1)^{(a+b+c)/2}\frac{[(a+b+c)/2+1]![(-a+b+c)/2]![(a-b+c)/2]![(a+b-c)/2]!}{[a]![b]![c]!}\cr
&&
\eea
The numbers in square brackets are known as $q$-number defined by
$$
[n] = \frac{q^{n/2}-q^{-n/2}}{q^{1/2}-q^{-1/2}}~,
$$
where $q\rightarrow 1$ gives the number $n$. Similarly, a $q$-factorial is given by 
$[n]! = [n][n-1] \ldots [3][2][1]$. 
\bigbreak
\subsection{Classical character variety}
As shown in Figure \ref{fig:theta2}, the tubular neighborhood of the theta graph is a genus-two handlebody. Hence, the moduli space  $\scM_{\rm flat}(\Sigma_{g=2},\SL(2,{\bC}))$ of $\SL(2,{\bC})$ flat connections on the boundary Riemann surface $\Sigma_{g=2}$ of genus-two is a hyper-K$\ddot{\rm a}$hler variety of $\dim_\bC=6$ \cite{Hitchin:1987} whose local coordinates  are spanned by the holonomy eigenvalues $(x_a, x_b, x_c)$ and twist coordinates $(y_a,y_b,y_c)$. 
The moduli space  $\scM_{\rm flat}(S^3\backslash \Theta,\SL(2,{\bC}))$ of $\SL(2,{\bC})$ flat connections on the complement $S^3\backslash \Theta$ is a Lagrangian subvariety of the moduli space $\scM_{\rm flat}(\Sigma_{g=2},\SL(2,{\bC}))$ with respect to the symplectic form $\omega=\frac{1}{\hbar}\sum_i d\log x_i\wedge d\log y_i$, which is described by the zero loci of  the set of algebraic equations $A_j(\Theta;y_j,\{x_i\})=0$ $(j=a,b,c)$
\bea\label{moduli-two}
 &&\scM_{\rm flat}(S^3\backslash \Theta,\SL(2,{\bC}))
 \cr&=& \Big\{(x_a,x_b,x_c,y_a,y_b,y_c)\in \scM_{\rm flat}(\Sigma_{g=2},\SL(2,{\bC}))\Big|A_j(\Theta;y_j,\{x_i\})=0 \ \ (j=a,b,c)\Big\}. 
\eea

The generalized volume conjecture states that the set of equations can be obtained by taking large color limit of the $U_q(\fraksl_2)$ quantum invariants of the theta graph.
In order to explicitly see the large color behavior, we use the asymptotic form of $[n]!$  in the limit of $q=e^\hbar\to 1$, $n\to\infty$ with $q^n=x_n$ fixed
\bea
[n]!\overset{\hbar\rightarrow 0,n\rightarrow\infty}{\longrightarrow}e^{\frac{1}{\hbar}g(x_n)}~,\quad {\rm where}\quad g(x_n) =-\frac{(\log{x_n})^2}{4}-\Li_2(x_n)~~. 
\label{gx} 
\eea
Thus, in the limit
\be
 \hbar\to 0,\quad (a,b,c)\to\infty,\quad x_a=q^{a/2}={\rm{fixed}},\quad  x_b=q^{b/2}={\rm{fixed}},\quad  x_c=q^{c/2}={\rm{fixed}},
 \label{limit1}
\ee
the quantum invariant of the theta graph \eqref{theta-inv} takes the form
\be
 J_{a,b,c}(\Theta) \sim e^{\tfrac{1}{\hbar}W(\Theta;x_a,x_b,x_c)},
 \label{asym-theta}
\ee
where the Neumann-Zagier potential is written by
\bea
W(\Theta;x_a,x_b,x_c)&=&\pi i \log{(x_a x_b x_c)}+g(x_a x_b x_c)+g(x_a^{-1} x_b x_c)+g(x_a x_b^{-1} x_c)+g(x_a x_b x_c^{-1})\cr
&&-g(x_a^2)-g(x_b^2)-g(x_  c^2)
\label{W-theta}
\eea
The set of equations which determines the moduli space $\scM_{\rm flat}(S^3\backslash \Theta,\SL(2,{\bC}))$ can be obtained from 
\be
y_j=\exp\left(x_j\frac{\partial W(\Theta;\{x_i\})}{\partial x_j}\right)~,
\ee
Hence, plugging the Neumann-Zagier potential \eqref{W-theta} into the equations, we conjecture that the moduli space $\scM_{\rm flat}(S^3\backslash \Theta,\SL(2,{\bC}))$ is determined by the zero loci of the following equations
\bea\label{A-two}
A_a(\Theta;y_a,\{x_i\})&=& x_a \left(x_a x_b-x_c\right) \left(x_a x_c-x_b\right) \left(1-x_a x_b x_c\right)+x_b x_c\left(1-x_a^2\right){}^2\left(x_a-x_b x_c\right)y_a~,\cr
A_b(\Theta;y_b,\{x_i\})&=&x_b \left(x_a x_b-x_c\right) \left(x_b x_c-x_a\right) \left(1-x_a x_b x_c\right)+x_a  x_c \left(1-x_b^2\right)^2\left(x_b-x_a x_c\right) y_b~,\cr
A_c(\Theta;y_c,\{x_i\})&=&x_c \left(x_a x_c-x_b\right) \left(x_b x_c-x_a\right) \left(1-x_a x_b x_c\right)+x_a x_b \left(1-x_c^2\right){}^2\left(x_c-x_a x_b\right)y_c~,\cr
&&
\label{A-poly-theta}
\eea
in the moduli space $\scM_{\rm flat}(\Sigma_{g=2},\SL(2,{\bC}))$. In fact, we have checked that the subvariety defined by the zero loci of the set of the equations \eqref{A-two} satisfy the Lagrangian condition. Namely, the symplectic form $\omega=\frac{1}{\hbar}\sum_i d\log x_i\wedge d\log y_i$ vanishes on the subvariety. We conjecture that the set \eqref{A-poly-theta} of the equations  locally determines the moduli space $\scM_{\rm flat}(S^3\backslash \Theta,\SL(2,{\bC}))\subset \scM_{\rm flat}(\Sigma_{g=2},\SL(2,{\bC}))$.
It is desirable to obtain the $\SL(2,\bC)$ character variety of the fundamental group $\pi_1(S^3\backslash \Theta)$ in terms of the complex Fenchel-Nielsen coordinates $\{x_i,y_i\}$ by the other means.

\subsection{Quantum character variety}
In order to quantize the classical character varieties $A_i(\{x_i\},y_i)$, we introduce 
operators $\hat{x}_i$ and $\hat{y}_i$, whose actions on $J_{a,b,c}(\Theta;q)$ are as
follows:
\bea
\hat{x}_a J_{a,b,c}(\Theta,q)&=&q^{a/2} J_{a,b,c}(\Theta,q)\cr
\hat{y}_a J_{a,b,c}(\Theta,q)&=& J_{a+2,b,c}(\Theta,q)
\eea
From the following recursion relation for $J_{a,b,c}(\Theta,q)$ 
\be
 J_{a+2,b,c}(\Theta,q)=-\frac{[(a+b+c)/2+2][(a-b+c)/2+1][(a+b-c)/2+1]}{[(-a+b+c)/2][a+1][a+2]} J_{a,b,c}(\Theta,q)
\label{thetarecur}~, 
\ee
we can deduce the form for the quantum character variety as 
\bea
\wh{A}_a(\Theta;\hat{y}_a,\{\hat{x}_i\})&=&\hat{x}_a \left(q \hat{x}_a \hat{x}_b-\hat{x}_c\right) \left(q \hat{x}_a \hat{x}_c-\hat{x}_b\right) \left(1-q^2 \hat{x}_a \hat{x}_b \hat{x}_c\right)\cr
&&+q^{-3/2}\hat{x}_b \hat{x}_c \left(1-\hat{x}_a^2\right) \left(q-\hat{x}_a^2\right) \left(\hat{x}_a-q \hat{x}_b \hat{x}_c\right) \hat{y}_a~,\cr
\wh{A}_b(\Theta;\hat{y}_b,\{\hat{x}_i\})&=&\hat{x}_b \left(q \hat{x}_a \hat{x}_b-\hat{x}_c\right) \left(q \hat{x}_b \hat{x}_c-\hat{x}_a\right) \left(1-q^2 \hat{x}_a \hat{x}_b \hat{x}_c\right)\cr
&&+q^{-3/2}\hat{x}_a \hat{x}_c \left(1-\hat{x}_b^2\right) \left(q-\hat{x}_b^2\right) \left(\hat{x}_b-q \hat{x}_a \hat{x}_c\right) \hat{y}_b~,\cr
\wh{A}_c(\Theta;\hat{y}_c,\{\hat{x}_i\})&=&\hat{x}_c \left(q \hat{x}_a \hat{x}_c-\hat{x}_b\right) \left(q \hat{x}_b \hat{x}_c-\hat{x}_a\right) \left(1-q^2 \hat{x}_a \hat{x}_b \hat{x}_c\right)\cr
&&+q^{-3/2}\left(\hat{x}_a \hat{x}_b 1-\hat{x}_c^2\right) \left(q-\hat{x}_c^2\right) \left(\hat{x}_c-q \hat{x}_a \hat{x}_b\right) \hat{y}_c~.
\label{qapoly}
 \eea
Our proposal of the quantization of the character variety for the theta graph is evident from the equality of 
$q\rightarrow 1$ limit of \eqref{qapoly} with the classical character varieties \eqref{A-poly-theta}. 

\section{Tetrahedron graph}
The colored quantum invariants  $J_{j_1,\ldots j_4,j_{12},j_{23}} (\smalltet{.5}{0.15};q)$ for the tetrahedron graph are know as the quantum $6j$-symbols for $U_q(\fraksl_2)$  \cite{Kirillov:1989}, up to a certain factor, which are written as state sum of product of four $q$-CG coefficients. In this paper, we adopt the notation used in \cite{masbaum1994}
\begin{figure}[h]
\begin{center}
\includegraphics[scale=.2]{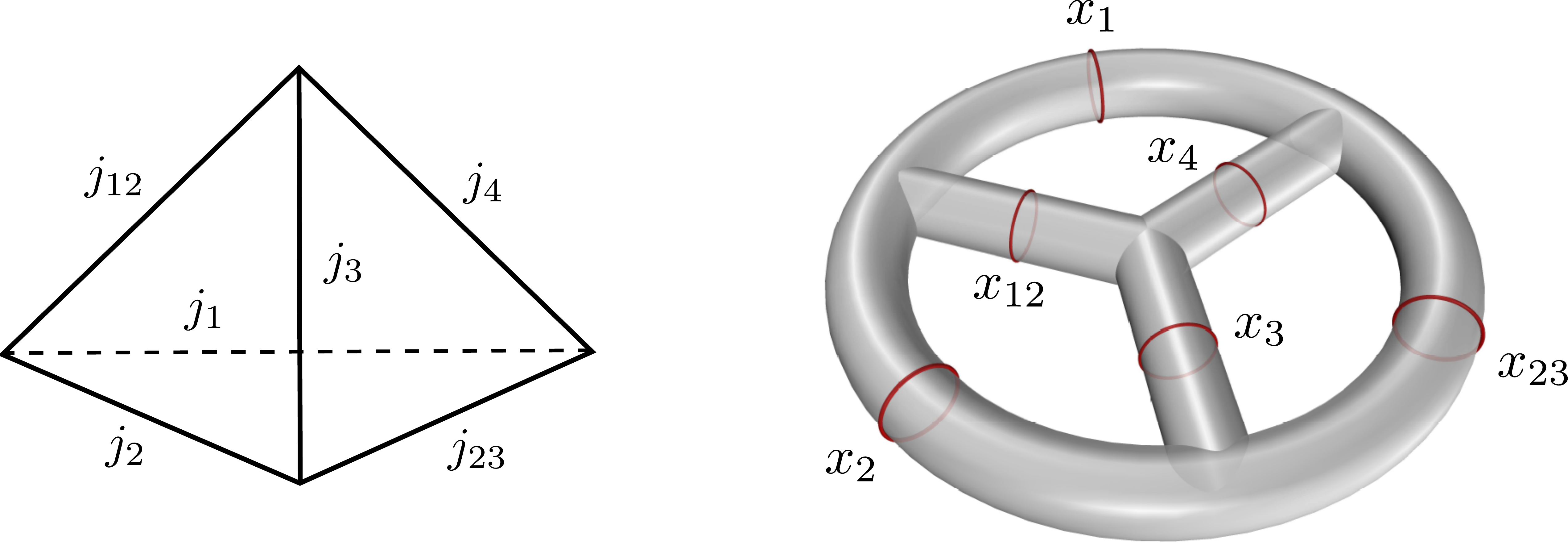}
\end{center}
\caption{Tetrahedron graph and its tubular neighborhood.}
\label{fig:tetrahedron}
\end{figure}
\bea
\left< \begin{matrix} j_1 & j_2 & j_{12} \\ j_3 & j_4 & j_{23}\end{matrix} \right> &=& \frac{\Delta
(j_1,j_2, j_{12}) \Delta(j_3,j_4,j_{12}) \Delta(j_1,j_4,j_{23}) \Delta
(j_2,j_3,j_{23})}{[j_1]![j_2]![j_3]![j_4]![j_{12}]![j_{23}]!}\cr
&& \times \sum_{m \ge 0} (-)^m [m+1]! ~~ {\bf \{ } [m-\tfrac12(j_1 - j_2 -j_{12})]! \cr
&& \times [m -\tfrac12 (j_3 - j_4 - j_{12})]! [m -\tfrac12( j_1 - j_4 - j_{23})]!\nonumber  \\
&& \times [m - \tfrac12(j_2 - j_3 - j_{23})]! [\tfrac12(j_1 + j_2 + j_3 + j_4) - m]!
\cr
&&\times  [\tfrac12(j_1+j_3+j_{12}+j_{23})-m]! [\tfrac12(j_2+j_4+j_{12}+j_{23})-m]!
 { \bf \} } ^{-1}
 \label{6j}
\eea
where we define
\be
\Delta(a,b,c) = [(-a+b+c)/2]! [(a-b+c)/2]! [(a+b-c)/2]!~.
\ee
Note that  $j_i \in \bZ$ represent the spin-$j_i/2$ representations of $U_q(\fraksl_2)$ assigned to the six edges (Figure \ref{fig:tetrahedron}). If  a singlet representation is placed on an edge of the tetrahedron, the tetrahedron invariant reduces to the quantum invariant \eqref{theta-inv} of the theta graph.
It is well-known that the tetrahedron invariant exhibits the following symmetry properties:
\be 
\left< \begin{matrix} j_1 & j_2 & j_{12} \\ j_3 & j_4 & j_{23}\end{matrix} \right> = \left< \begin{matrix} j_2 & j_1 & j_{12} \\ j_4 & j_3 & j_{23}\end{matrix} \right> = \left< \begin{matrix} j_{12} & j_2 & j_1 \\ j_{23} & j_4 & j_3\end{matrix} \right> = \left< \begin{matrix} j_3 & j_2 & j_{23} \\ j_1 & j_4 & j_{12}\end{matrix} \right> = \left< \begin{matrix} j_3 & j_4 & j_{12} \\ j_1 & j_2 & j_{23}\end{matrix} \right>~~. 
\label{symmetry}
\ee
In what follows, we shall obtain the classical and quantum character varieties of the complement of the tetrahedron graph in $S^3$. To this end, we will ignore the factor in front of the summation in \eqref{6j} since the effect of the factor can be absorbed into the redefinition of the variables $y_i$. It turns out  \cite{Kirillov:1989} that the summation part can be written in terms of the basic hypergeometric series
\bea
&&\left< \begin{matrix} j_1 & j_2 & j_{12} \\ j_3 & j_4 & j_{23}\end{matrix} \right>^\prime  \cr
&=&\tfrac{(-)^{(j_1+j_2+j_3+j_4)/2}[\half(j_1+j_2+j_3+j_4)+1]!}{[\half(j_1 + j_2 - j_{12})]![\half(j_1 + j_4 - j_{23})]![\half(j_3 + j_4 - j_{12})]![\half(j_2 + j_3 -j_{12})]![\half(j_{12}+j_{23}-j_2-j_4)]![\half(j_{12}+j_{23}-j_1-j_2)]!}\cr
 &&{}_4\varphi_3\left(\begin{array}{c}q^{\half(-j_1 - j_2 + j_{12})},~ q^{\half(-j_1 - j_4 + j_{23})},~q^{\half(-j_3 - j_4 + j_{12})},~q^{\half(-j_2 - j_3 + j_{12})} \\ q^{-\half(j_1+j_2+j_3+j_4)-1},~q^{\half(j_{12}+j_{23}-j_2-j_4)+1},~q^{\half(j_{12}+j_{23}-j_1-j_2)+1} \end{array} \Big|q,q\right)
 \label{6j-2}
\eea

\subsection{Classical character variety}
It is easy to see that the tubular neighborhood of the tetrahedron graph is a genus-3 handlebody. Hence, as Figure \ref{fig:tetrahedron} depicted, the moduli space $\scM_{\rm flat}(\Sigma_{g=3},\SL(2,{\bC}))$ of $\SL(2,\bC)$ flat connections on the boundary Riemann surface $\Sigma_{g=3}$ will be a 12-dimensional hyper-K$\ddot{\rm a}$hler variety \cite{Hitchin:1987} whose coordinates are spanned by the holonomy eigenvalues $x_i \ (i\in\{1,\cdots,4,12,23\})$ and their conjugate variables $y_i\ (i\in\{1,\cdots,4,12,23\})$. The classical character variety $\scM_{\rm flat}(S^3\backslash\smalltet{.4}{0.09},\SL(2,{\bC}))$  of the complement of the tetrahedron graph is determined by a set of six Laurent polynomials 
\be\label{classical-theta}
 A_k (\smalltet{.5}{0.15}; y_k, \{x_i\}) = 0 ~ , \quad  k\in\{1,\cdots,4,12,23\}~,
\ee
which cuts out the Lagrangian subvariety of the moduli space $\scM_{\rm flat}(\Sigma_{g=3},\SL(2,{\bC}))$. To see the explicit form of \eqref{classical-theta} via the generalized volume conjecture, we will take the limit 
\be
q=e^{\hbar}\to 1\,, \qquad j_i\to\infty\,, \qquad x_{i}=q^{j_i/2}=\rm{fixed}~,\, \quad i\in\{1,\cdots,4,12,23\}
\label{limit2}
\ee
of the colored quantum invariants of the tetrahedron graph. In this limit, the summation over  $m$  in \eqref{6j-2} is approximated an integral over a variable $z=q^{m}$ 
\be
\left< \begin{matrix} j_1 & j_2 & j_{12} \\ j_3 & j_4 & j_{23}\end{matrix} \right>^\prime\sim \int e^{\frac{1}{\hbar}\widetilde\cW(\smalltet{.4}{0.10};\{x_i\};z)}dz,
\label{asym-6jnorm}
\ee
where the Neumann-Zagier potential can be read off from \eqref{gx} as
\bea
\widetilde\cW(\smalltet{.5}{0.15};\{x_i\},z)&=&\pi i \log(z)+g(z)\cr
&& -g\left(\tfrac{z}{x_{j_1}x_{j_2}x_{j_{12}}}\right)-g\left(\tfrac{z}{x_{j_3}x_{j_4}x_{j_{12}}}\right)-g\left(\tfrac{z}{x_{j_1}x_{j_4}x_{j_{23}}}\right)-g\left(\tfrac{z}{x_{j_2}x_{j_3}x_{j_{23}}}\right)\cr
&&-g\left(\tfrac{x_{j_1}x_{j_2}x_{j_3}x_{j_4}}{z}\right)-g\left(\tfrac{x_{j_1}x_{j_{12}}x_{j_3}x_{j_{23}}}{z}\right)-g\left(\tfrac{x_{j_2}x_{j_{12}}x_{j_4}x_{j_{23}}}{z}\right)
\label{superpotential}
\eea
The leading contribution of the asymptotics comes from the saddle point $z=z_0$ 
\be
\left.\frac{\partial \widetilde\cW(\smalltet{.5}{0.15};\{x_i\},z)}{\partial z}\right|_{z=z_0}=0~.
\label{supercurve1}\\
\ee
Therefore, the zero locus of the classical character variety is determined by
\be
y_k=\exp\left(x_k\frac{\partial \widetilde\cW(\smalltet{.5}{0.15};\{x_i\},z_0)}{\partial x_k}\right)~. 
\ee
Substituting the Neumann-Zagier potential \eqref{superpotential} into the above equations for $k=1$, we have
\begin{footnotesize}
\bea\label{saddle_point} 
1&=&\frac{(z_0-1)(z_0-x_{1}x_{2}x_{3}x_{4})(z_0-x_{1}x_{3} x_{12}x_{23})(z_0-x_{2} x_{4} x_{12} x_{23})}{(z_0-x_{1} x_{2} x_{12})(z_0-x_{1} x_{4} x_{23})(z_0-x_{2} x_{3} x_{23})(z_0-x_{3} x_{4} x_{12})}~,\cr
y_1&=&\frac{x_3 \left(z_0-x_1 x_4 x_{23}\right) \left(z_0-x_1 x_2 x_{12}\right)}{\left(z_0-x_1 x_2 x_3 x_4\right) \left(z_0-x_1 x_3 x_{12} x_{23}\right)}~. 
\eea
\end{footnotesize}
Eliminating $z_0$ from these two equations, we obtain an explicit expression
\begin{small}
\bea\label{cCV}
A_1(\smalltet{.5}{0.15};y_1,\{x_i\})
 &=&x_3 \left(x_1 x_4-x_{23}\right) \left(x_4-x_1 x_{23}\right) \left(x_1 x_2-x_{12}\right) \left(x_2-x_1 x_{12}\right)y_1^2\cr
 &&-\Big[x_4x_{12}(1-x_1^2)^2  (x_3-x_2 x_{23})(1-x_3 x_2 x_{23}) \cr
 &&+x_3(x_{23}-x_1 x_4)(x_1-x_4 x_{23})(x_1-x_2 x_{12})(x_2-x_1 x_{12})\cr
 &&+x_3 (x_4-x_1 x_{23})(1-x_1 x_4 x_{23})(x_{12}-x_1 x_2)(1-x_1 x_2 x_{12})\Big]y_1\cr
 &&+x_3 \left(x_1-x_4 x_{23}\right) \left(x_1-x_2 x_{12}\right) \left(1-x_1 x_4 x_{23}\right) \left(1-x_1 x_2 x_{12}\right)~.
\eea
\end{small}
From this expression, we can find the other formulas by using the symmetry property \eqref{symmetry}, 
\bea
A_2(\smalltet{.5}{0.15};y_2,\{x_i\})&=&A_2(\smalltet{.5}{0.15};y_1\rightarrow y_2,x_2\leftrightarrow x_1,x_3\leftrightarrow x_4)~,\cr
A_3(\smalltet{.5}{0.15};y_3,\{x_i\})&=&A_1(\smalltet{.5}{0.15};y_1\rightarrow y_3,x_3\leftrightarrow x_1,x_{23}\leftrightarrow x_{12})~,\cr
A_4(\smalltet{.5}{0.15};y_4,\{x_i\})&=&A_1(\smalltet{.5}{0.15};y_1\rightarrow y_4,x_3\leftrightarrow x_2,x_1\leftrightarrow x_4)~,\cr
A_{12}(\smalltet{.5}{0.15};y_{12},\{x_i\})&=&A_1(\smalltet{.5}{0.15};y_1\rightarrow y_{12},x_3\leftrightarrow x_{23},x_1\leftrightarrow x_{12})~,\cr
A_{23}(\smalltet{.5}{0.15};y_{23},\{x_i\})&=&A_1(\smalltet{.5}{0.15};y_1\rightarrow y_{23},x_3\leftrightarrow x_{12},x_1\leftrightarrow x_{23})~.
\eea
Thus, we conjecture that the zero locus of these six polynomial equations represents the classical character variety of the tetrahedron graph in the moduli space $\scM_{\rm flat}(\Sigma_{g=3},\SL(2,{\bC}))$. Indeed, we have verified that it is subject to the Lagrangian condition by showing the symplectic form
$\omega={i\over \hbar}\sum_i d \log x_i \wedge d\log y_i=0$ on the 
subvariety. 

\subsection{Quantum character variety }
As we have seen, the colored quantum invariants of the tetrahedron graph can be written in terms of a balanced basic hypergeometric series  ${}_4\varphi_3$. Moreover, it is a certain orthogonal polynomial called an Askey-Wilson polynomial \cite{Askey:1979}. Thus, we can obtain the recursion relation \cite{Askey:1979,kachurik} which is explicitly  written as
\be
\alpha\left < \begin{matrix} j_1+2 & j_2 & j_{12} \\ j_3 & j_4 & j_{23}\end{matrix} \right >^\prime-\beta\left < \begin{matrix} j_1 & j_2 & j_{12} \\ j_3 & j_4 & j_{23}\end{matrix} \right >^\prime+\gamma\left< \begin{matrix} j_1-2 & j_2 & j_{12} \\ j_3 & j_4 & j_{23}\end{matrix} \right >^\prime=0~,
\label{recursum}
\ee
where
\bea
\alpha&=&[j_1][(j_1+j_4-j_{23})/2+1][(j_1+j_{23}-j_4)/2+1][(j_1+j_2-j_{12})/2+1][(j_1+j_{12}-j_2)/2+1]\cr 
\beta&=&[j_1+2][(j_1+j_{23}-j_4)/2][(j_1+j_4+j_{23})/2+1][(j_1+j_2-j_{12})/2][(j_1+j_2+j_{12})/2+1]\cr
&&~~~~+[j_1][(j_1+j_4-j_{23})/2+1][(j_{23}+j_4-j_1)/2][(j_2+j_{12}-j_1)/2][(j_1+j_{12}-j_2)/2+1]\cr
&&~~~~-[j_1][j_1+1][j_1+2][(j_3+j_2+j_{23})/2+1][(j_2+j_{23}-j_3)/2]\cr
\gamma&=&[j_1+2][(j_1+j_4+j_{23})/2+1][(j_{23}+j_4-j_1)/2+1][(j_2+j_{12}+j_1)/2+1]\cr
&&\times[(j_2+j_{12}-j_1)/2+1].
\eea


Using the following action of the quantum operator $\hat x_i$ and the conjugate operator $\hat y_1$ for the tetrahedron invariant 
\bea
\hat{x}_i \left\langle \begin{matrix} j_1 & j_2 & j_{12} \\ j_3 & j_4 & j_{23}\end{matrix} \right\rangle^\prime&=&q^{\frac{j_i}{2}}\left\langle \begin{matrix} j_1 & j_2 & j_{12} \\ j_3 & j_4 & j_{23}\end{matrix} \right\rangle^\prime\rm\\
\hat{y}_1 \left\langle \begin{matrix} j_1 & j_2 & j_{12} \\ j_3 & j_4 & j_{23}\end{matrix} \right\rangle^\prime&=&\left\langle \begin{matrix} j_1+2 & j_2 & j_{12} \\ j_3 & j_4 & j_{23}\end{matrix} \right\rangle^\prime\rm,
\eea
we can write \eqref{recursum} as an operator which annihilates the quantum invariants of the tetrahedron graph
\bea
\label{qCV}
&&\wh{A}_1(\smalltet{.5}{0.15};\hat{y}_1,\{\hat{x}_i\};q)\cr&=&
q^3 \hat{x}_3 \left(q^2-\hat{x}_1^2\right)  \left(\hat{x}_1 \hat{x}_4-\hat{x}_{23}\right) \left(\hat{x}_4-\hat{x}_1 \hat{x}_{23}\right) \left(\hat{x}_1 \hat{x}_2-\hat{x}_{12}\right) \left(\hat{x}_2-\hat{x}_1 \hat{x}_{12}\right)\hat{y}_1^2\cr
&&-\Big[q^2\hat{x}_4\hat{x}_{12}(1-\hat{x}_1^2)(1-q \hat{x}_1^2)(1-q^2 \hat{x}_1^2) (\hat{x}_3-\hat{x}_2 \hat{x}_{23})(1-q \hat{x}_3 \hat{x}_2 \hat{x}_{23}) \cr
&&+q^3 \hat{x}_3(1-\hat{x}_1^2)(\hat{x}_{23}-q \hat{x}_1 \hat{x}_4)(\hat{x}_1-\hat{x}_4 \hat{x}_{23}) (\hat{x}_1-\hat{x}_2 \hat{x}_{12})(\hat{x}_2-q \hat{x}_1 \hat{x}_{12})\cr
&&+q\hat{x}_3(1-q^2 \hat{x}_1^2)(\hat{x}_4-\hat{x}_1 \hat{x}_{23})(1-q \hat{x}_1 \hat{x}_4 \hat{x}_{23})(\hat{x}_{12}-\hat{x}_1 \hat{x}_2)(1-q \hat{x}_1 \hat{x}_2 \hat{x}_{12})\Big]\hat{y}_1\cr
&&+\hat{x}_3 \left(1-q^4 \hat{x}_1^2\right)  \left(\hat{x}_1-\hat{x}_4 \hat{x}_{23}\right) \left(1-q^2 \hat{x}_1 \hat{x}_4 \hat{x}_{23}\right) \left(\hat{x}_1-\hat{x}_2 \hat{x}_{12}\right)\left(1-q^2 \hat{x}_1 \hat{x}_2 \hat{x}_{12}\right)~.
\eea
As $q\to 1$, the equation \eqref{qCV} reduces to its classical counterpart \eqref{cCV} up to a factor $\left(1-x_1^2\right)$. The other five equations easily follows from the symmetry property \eqref{symmetry}
\bea  
\wh{A}_2(\smalltet{.5}{0.15};\hat{y}_2,\{\hat{x}_i\},q)&=&\wh{A}_1(\smalltet{.5}{0.15};\hat{y}_1\rightarrow \hat{y}_2,\hat{x}_2\leftrightarrow \hat{x}_1,\hat{x}_3\leftrightarrow \hat{x}_4)~,\cr
\wh{A}_3(\smalltet{.5}{0.15};\hat{y}_3,\{\hat{x}_i\},q)&=&\wh{A}_1(\smalltet{.5}{0.15};\hat{y}_1\rightarrow \hat{y}_3,\hat{x}_3\leftrightarrow \hat{x}_1,\hat{x}_{23}\leftrightarrow \hat{x}_{12})~,\cr
\wh{A}_4(\smalltet{.5}{0.15};\hat{y}_4,\{\hat{x}_i\},q)&=&\wh{A}_1(\smalltet{.5}{0.15};\hat{y}_1\rightarrow \hat{y}_4,\hat{x}_3\leftrightarrow \hat{x}_2,\hat{x}_1\leftrightarrow \hat{x}_4)~,\cr
\wh{A}_{12}(\smalltet{.5}{0.15};\hat{y}_{12},\{\hat{x}_i\},q)&=&\wh{A}_1(\smalltet{.5}{0.15};\hat{y}_1\rightarrow \hat{y}_{12},\hat{x}_3\leftrightarrow \hat{x}_{23},\hat{x}_1\leftrightarrow \hat{x}_{12})~,\cr
\wh{A}_{23}(\smalltet{.5}{0.15};\hat{y}_{23},\{\hat{x}_i\},q)&=&\wh{A}_1(\smalltet{.5}{0.15};\hat{y}_1\rightarrow \hat{y}_{23},\hat{x}_3\leftrightarrow \hat{x}_{12},\hat{x}_1\leftrightarrow \hat{x}_{23})~.
\eea
In fact, these recurrence relations can be considered as the actions of line operators on the $S$-duality wall of $\cN=2$ $SU(2)$ supersymmetric gauge theories with $N_f=4$ \cite{Dimofte:2013lba}.

\section{Discussions}
In this letter, we have proposed extensions of two conjectures: (i) the generalized volume conjecture for trivalent graphs (ii) the AJ conjecture for trivalent graphs.
Particularly, our elaborate exercise for the theta and tetrahedron graphs provides a strong evidence for both the conjectures. To confirm that both the conjectures are true in these cases, it requires to derive the expressions of the classical character variety  by gluing equations of ideal tetrahedra triangulations of the complement. Even if the conjectures are true, the configurations of character varieties obtained in this paper provide only local information. The further investigation has to be undertaken to study their global structure.

Although we have dealt with only the theta and tetrahedron graphs in this paper, one can also consider general planar trivalent graphs as well as  knotted trivalent graphs \cite{Roland}. However, we do \emph{not} expect that the large color asymptotics and recursion relations with respect to only one color would provide character varieties of any trivalent graphs. In fact, in the context of multi-component links, recursion relations involving multi-colors have to be taken into account in order to determine the character variety of the complement in general. Therefore, we rather expect that there are $E$ recursion relations involving multi-colors 
\bea
\wh A_i(\Gamma;\{\hat x_j\},\{\hat y_k\};q)  J_{n_1,\cdots,n_E}(\Gamma;q) =0~,\qquad (i=1,\cdots,E).
\eea
whose classical limits determine the character variety of a trivalent graph. 

In the case of knots, both the conjectures are extended for the HOMFLY polynomials colored by symmetric representations \cite{Aganagic:2012jb,Fuji:2012nx,Garoufalidis:2012rt}. To the contrary, in the case of $U_q(\fraksl_N)$ quantum invariants of a trivalent graph, the fusion rule does not allow to increase just one of the colors for representations of $U_q(\fraksl_N)$. Hence, there is no recursion relation with respect to only one of the colors for higher rank quantum invariants of trivalent graphs. This is the reason why an expression for quantum $6j$-symbols of $U_q(\fraksl_N)$ \cite{Nawata:2013ppa} for a simplest class of multiplicity-free representations cannot be written as an Askey-Wilson polynomial. Consequently, the recursion relations for $U_q(\fraksl_N)$ quantum invariants of trivalent graphs are much more complicated involving multiple colors and general Young tableaux. The feasibility is beyond the scope of our current techniques.

It is well known that the quantum $6j$-symbols expressed in terms of the Askey-Wilson polynomials ${}_4\varphi_3$ exhibit remarkably rich symmetries and dualities. Translating this property in the context of the 3d/3d correspondence, the dualities turn out to be 3d mirror symmetries. According to \cite{Dimofte:2012pd}, there are 40 descriptions in terms of $\cN=2$ $U(1)$ gauge theory  and 32 descriptions in terms of $\cN=2$ $SU(2)$ gauge theory \cite{Teschner:2012em} for the 3d gauge theory corresponding to $6j$-symbols. The interpretation of these relationships in terms of Stokes phenomena of holomorphic blocks \cite{Beem:2012mb} is currently under investigation.

\section*{Acknowledgement}
S.N. would like to thank T. Dimofte, S. Gukov and R. van der Veen for valuable discussions. The work of S.N. is partially supported by the ERC Advanced Grant no.~246974, {\it``Supersymmetry: a window to non-perturbative physics''}.

\bibliographystyle{alpha}

\end{document}